\documentclass[a4paper,10pt]{article}
\usepackage[utf8]{inputenc}
\usepackage[utf8]{inputenc}
\usepackage[ngerman, english]{babel}
\usepackage{amsmath,amsthm,amssymb,amsfonts}
\usepackage{mathrsfs}
\usepackage{ifthen}
\usepackage{enumerate}
\usepackage{comment}
\usepackage{subcaption}
\usepackage{hyphsubst}
\usepackage{xcolor}
\usepackage{calc}
\usepackage{graphicx}
\usepackage[toc,page]{appendix}
\usepackage[mathscr]{euscript}

\makeatletter
\def\@hspace#1{\begingroup\setlength\dimen@{#1}\hskip\dimen@\endgroup}
\makeatother

\newtheorem{theorem}{Theorem}[section]
\newtheorem{lemma}[theorem]{Lemma}

\newtheorem{remark}[theorem]{Remark}

\renewcommand{\Im}{\operatorname{Im}}

\newcommand{\tv}{\boldtau}
\newcommand{\ol}[1]{\overline{#1}}

\newcommand{\spl}{\langle}
\newcommand{\spr}{\rangle}
\newcommand{\bpm}{\begin{pmatrix}}
\newcommand{\epm}{\end{pmatrix}}

\DeclareMathOperator{\curl}{curl}
\DeclareMathOperator{\Curl}{Curl}

\renewcommand{\div}{\operatorname{div}}

\newcommand{\boldtau}{\boldsymbol{\tau}}

\newcommand{\boldu}{\mathbf{u}}

\newcommand{\boldw}{\mathbf{w}}
\newcommand{\boldx}{\mathbf{x}}

\newcommand{\boldE}{\mathbf{E}}

\newcommand{\boldH}{\mathbf{H}}

\newcommand{\boldL}{\mathbf{L}}

\newcommand{\calE}{\mathcal{E}}

\newcommand{\calU}{\mathcal{U}}

\definecolor{brickred}{rgb}{0.8, 0.25, 0.33}
\definecolor{bostonuniversityred}{rgb}{0.8, 0.0, 0.0}
\definecolor{cornellred}{rgb}{0.7, 0.11, 0.11}
\definecolor{corn}{rgb}{0.98, 0.93, 0.36}
\definecolor{schoolbusyellow}{rgb}{1.0, 0.85, 0.0}
\definecolor{TUblue}{rgb}{0,102,153}
\colorlet{TUbluelight}{TUblue!30!white}

\usepackage[margin=2.5cm]{geometry}
\usepackage{url}
\def\bfx{\mathbf{x}}
\def\bfE{\mathbf{E}}
\def\bfH{\mathbf{H}}
\def\bfu{\mathbf{u}}
\def\bfU{\mathbf{U}}

\def\bfe{\mathbf{e}}
\def\bfnu{\boldsymbol{\nu}}
\def\bftv{\boldsymbol{\tau}}
\def\hboldE{\hat{\boldE}}
\def\hboldH{\hat{\boldH}}
\title{On the analysis of waveguide modes in an electromagnetic transmission line%
}
\author{Martin Halla%
\thanks{%
Institut f\"ur Numerische und Angewandte Mathematik,
Georg-Augst Universität Göttingen,
Lotzestr.\ 16-18, 37083 Göttingen, Deutschland.
e-mail: m.halla@math.uni-goettingen.de}
\and
Peter Monk%
\thanks{%
Department of Mathematical Sciences University of Delaware, Newark DE 19716, USA. e-mail: monk@udel.edu
}
}

\begin{document}

\maketitle
\abstract{Modal expansions are useful to understand wave propagation in an infinite electromagnetic  transmission line or waveguide.  They can also be used to construct generalized Dirichlet-to-Neumann maps that can be used to provide artificial boundary conditions for truncating  a computational domain when discretizing the field by finite elements.  The modes of a waveguide arise as eigenfunctions of a non-symmetric eigenvalue problem, and the eigenvalues determine the propagation (or decay) of the modes along the waveguide.  For the successful use of waveguide modes, it is necessary to know that the modes exist and form a dense set in a suitable function space containing the trace of the electric field in the waveguide.  This paper is devoted to proving such a density result using the methods of Keldysh.  We also show that the modes satisfy a useful orthogonality property, and show how the Dirichlet-to-Neumann map can be calculated.  Our existence and density results are proved under realistic regularity assumptions on the cross section of the waveguide, and the electromagnetic properties of the materials in the waveguide, so generalizing existing results.}

\section{Introduction}

\begin{figure}
\centering
\resizebox{0.45\textwidth}{!}{\includegraphics{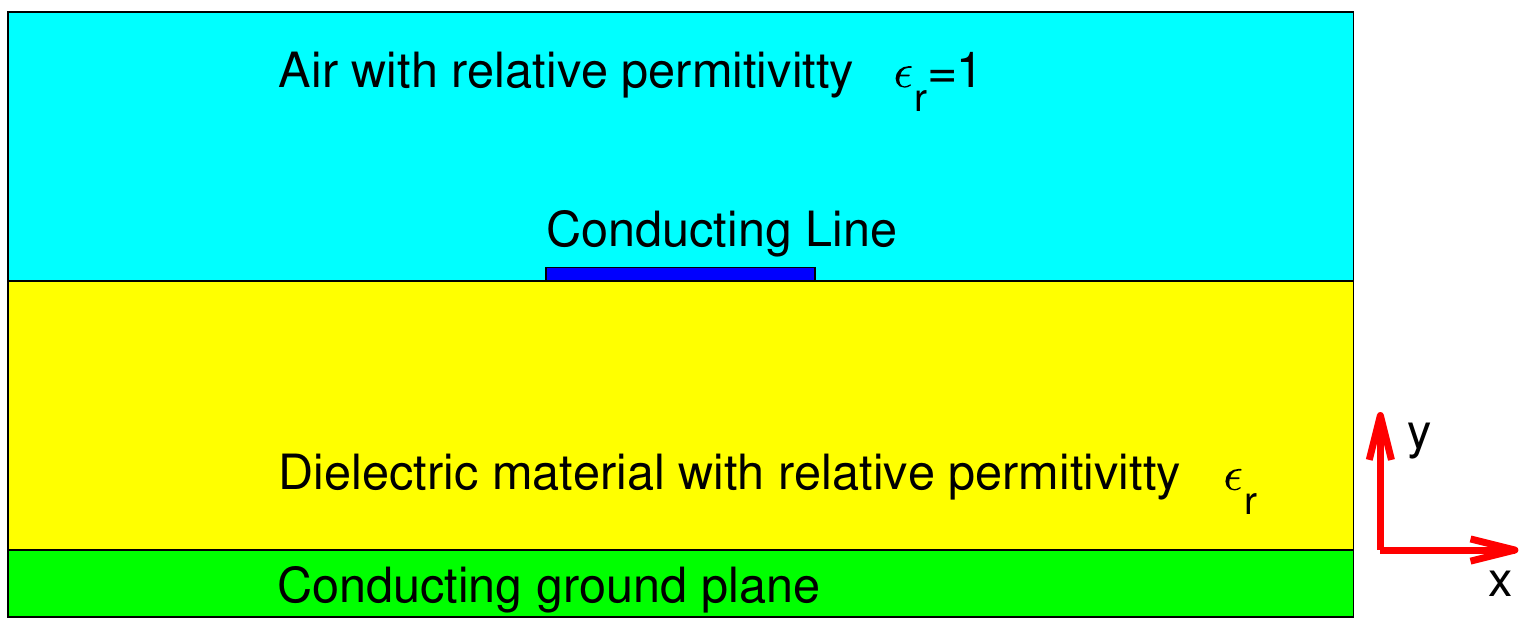}}\hspace{.2in}\resizebox{0.2\textwidth}{!}{\includegraphics{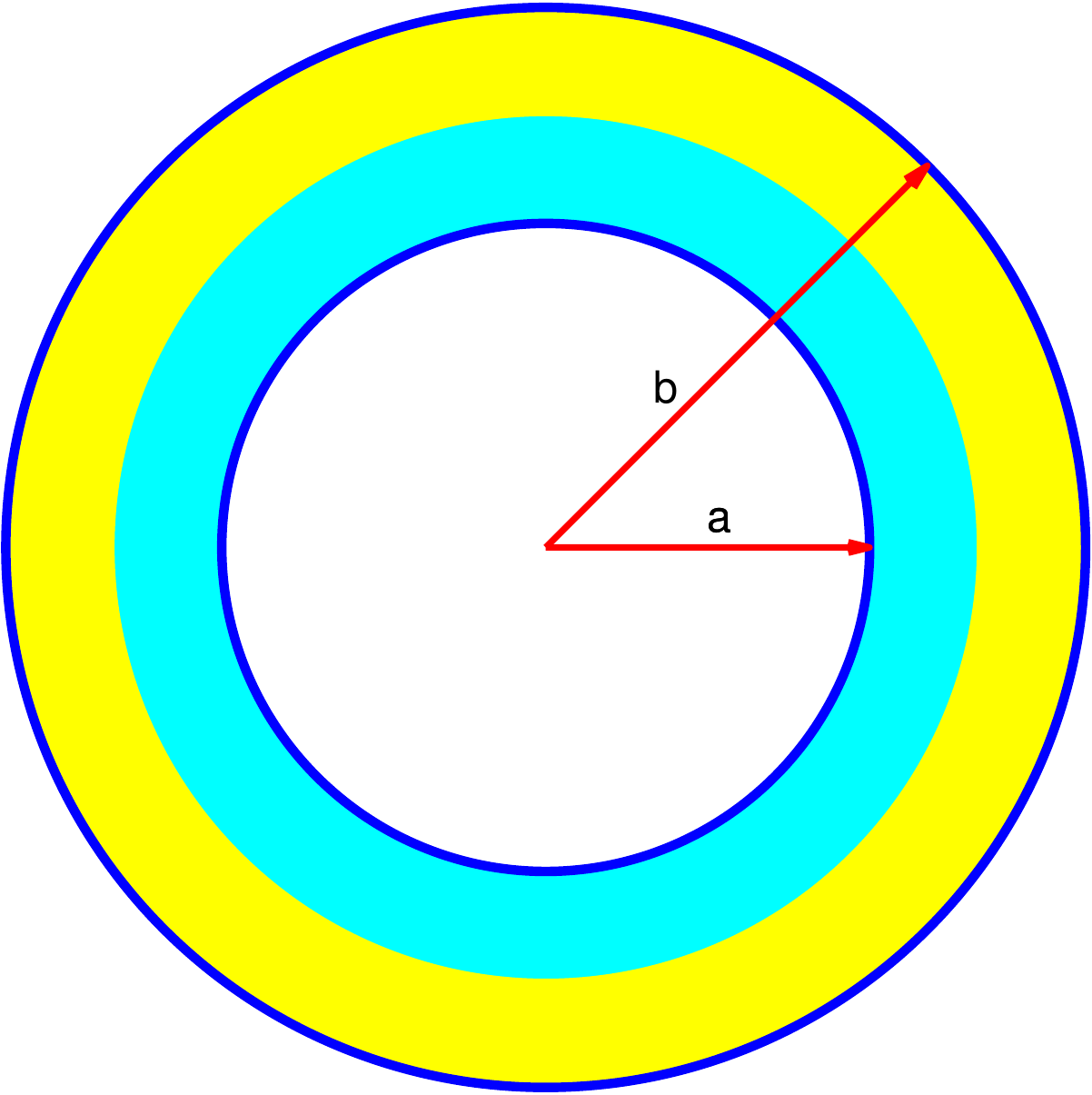}}
\caption{Left: Cartoon of a simple microstrip transmission line. The transmission line (dark blue region) is above a conducting ground plane (green). Right: Cartoon of a coaxial waveguide with a perfect conducting boundary (blue curves) and piecewise constant electromagnetic parameters indicated by the 
different shading.}
\label{Fig1}
\end{figure}

Transmission lines are devices used to transmit electromagnetic waves, for example connecting antenna elements in a device.  There are several different types of transmission line including coaxial waveguides, microstrips, and strip lines~\cite{Pozar:11}.  Out of the many possible designs, 
we shall focus on ones which have the following characteristics: they are long compared to the wavelength of the electromagnetic field with an invariant shape along the intended direction of propagation (assumed to be along the $z$-coordinate axis in this paper) and they consist of a conductor, or multiple conductors often over or surrounded by a Perfectly Electrically Conducting  (PEC) surface~\cite[Page 96]{Pozar:11}.  They may involve magnetic or dielectric components. The cross-section of two simple cartoon examples in the $(x,y)$-plane  is shown in Fig.~\ref{Fig1}.  In the left hand panel, a representative model of a microstrip transmission with width $W$, thickness $t$ and  positioned at height $h$ above a PEC ground plane is shown (see for example~\cite{Dassault,Lumerical,VardapetyanDemkowiczNeikirk:03}).  A dielectric substrate with 
electric permittivity $\epsilon$ separates the line from the ground plane. The right panel shows a cartoon of  the cross-section of a waveguide.

For the microstrip transmission line, the region above the ground plane is often infinite in extent (a half plane) and it is desired to compute modes that are ``localized'' near the transmission line.  Because such modes are assumed to decay rapidly away from the microstrip, it is usual in practice to truncate the computational structure
with a simple homogeneous PEC boundary condition~\cite{Dassault,Lumerical}.   After this truncation we see that in both cases in Fig.~\ref{Fig1} we have a translationally invariant domain bounded by a perfect conductor (for a discussion of the non-truncated problem, see \cite{JolyPoirier:99}).

For these structures the modes of the waveguide are associated to a propagation constant $\beta$ and corresponding angular frequency $\omega$ of the radiation (see the next section for a precise definition of these terms).   An early paper due to Bamberger and Bonnet~\cite{Bamberger:90} considered the problem of computing $\omega$ when $\beta$ is fixed and
where the optical fiber is modeled as an infinitely long cylindrical waveguide, with variable permittivity in an infinite domain. The presence of variable coefficients and perfect conductors in microstrips may produce solutions to Maxwell's equations with low regularity.
The study of guided waves is continued in \cite{JolyPoirier:95} where again  $\omega$ is computed as a function of $\beta$
and the resulting dispersion relation is studied. Computation of the dispersion relations using a non-conforming method is studied in \cite{JolyPoirierRobertsTrouve:96}.

Our paper differs from the above in that we assume $\omega$ to be fixed and  $\beta$ to be the eigenvalue. We also allow the presence of thin PEC structures (producing fields with corner singularities)
as well as allowing both electric permittivity and magnetic permeability to vary.  
Our goal is then to prove the density of the eigen-modes in an appropriate space. To prove this density result, we apply \cite[Theorem~2.1]{Halla:22MSteklovStab} (which in turn traces back to the work of Keldysh and is a reformulation of \cite[Chapter V, Theorem 8.3]{GohbergKrein:69} and \cite[Theorem 4.2, Corollary 3.3]{Markus:88}) to a variational formulation of the modal eigenvalue problem  suggested by Vardepetyan and Demkowicz~\cite{VardapetyanDemkowicz:03}.  

Such a density result is needed to give a theoretical underpinning to modal methods for solving
waveguide problems.
In addition, if we wish to compute the solution of a scattering problem in the waveguide in which there is a locally perturbed geometry or structure (for example a change in diameter of the waveguide, a corner or junction or a flaw) it is necessary to truncate the domain to allow finite elements to be used in a bounded region containing the perturbation, while a suitable truncation condition is used to account for the infinite translationally invariant parts of the structure.  One way to do this is to use the electromagnetic version of the Dirichlet-to-Neumann map on a cross section of the waveguide.  For this to be successful we need to know that, for a fixed $\omega$, the set of modes corresponding to all propagation constants are dense in the appropriate boundary space.  

In \cite{VardapetyanDemkowicz:03} it is shown both theoretically and computationally that the method proposed there gives a convergent approximation of the
the propagation constants. Numerical tests and hp-implementation of this method can be found in \cite{VardapetyanDemkowiczNeikirk:03}. An alternative approach is to eliminate the four fields transverse to the direction of the
transmission line.  This has been studied for the  problem of computing modes in a waveguide with piecewise constant electromagnetic properties~\cite{SheSmi13,SheSmi13a}.  A representative geometry is shown in Fig.~\ref{Fig1} right panel. Assuming the electromagnetic coefficients are piecewise constant, an eigenvalue problem for a 4th order operator pencil can be derived.  It is proved that eigenvalues exist~\cite{SheSmi13} and that only finitely many are real, and in~\cite{SheSmi13a} completeness properties are also derived.  It seems to be critical for this approach that the coefficients are piecewise constant, whereas the approach in \cite{VardapetyanDemkowicz:03} allows for general
piecewise smooth (precisely piecewise $W^{1,\infty}$) coefficients, so that our results extend the completeness theory of~\cite{SheSmi13,SheSmi13a}
to allow for graded materials.

The remainder of the paper proceeds as follows. In Section ~\ref{deriv} we start by a formal derivation of  Vardapetyan and Demkowicz's eigenvalue problem from Maxwell's equations.  In particular we state the version that we shall analyze.  We then proceed to  prove our main theorem, the existence and density result, in Section~\ref{density}.  Then in Section~\ref{DtN} we show how the modes can be used to compute a DtN operator suitable for discretization.

Concerning notation,  we shall denote by $\Omega\subset\mathbb{R}^2$ the cross section of the transmission line enclosed by a bounded PEC surface. Additional assumptions on $\Omega$ will be formulated in the theorems we state.
We denote by $\bfnu$ the two dimensional outward normal to $\Omega$ and by $\bftv$  the unit tangential vector to the domain on boundary $\partial\Omega$.
For square integrable scalar or vector functions $f$ and $g$ we use the notation
\[
\spl f,g\spr=\int_{\Omega}f\overline{g}\,dx.
\]
The same notation is also used for duality pairings.

\section{Derivation of the modal eigenvalue problem}\label{deriv}
We now provide a brief derivation of the eigenvalue problem due to Vardepetyan and Demkowicz~\cite{VardapetyanDemkowicz:03}.
The axis of the waveguide is assumed to be parallel to the $z$ coordinate axis and we let $\Omega\subset \mathbb{R}^2$ denote the cross-section of the waveguide in $(x,y)$ plane.  Critically the magnetic permeability $\mu$ and the electric permeability $\epsilon$ of the material in the waveguide are assumed to be independent of $z$. Usually $\mu$ is constant, but in case of magnetic components in the device it may differ from unity.  For simplicity, we shall assume that $\mu$ is a scalar quantity. Both  $\epsilon$ and $\mu$ are assumed to be  real valued, bounded and uniformly positive scalar functions of position (precisely they are assumed to be piecewise $W^{1,\infty}$ functions).

Let $\hat\bfx=(x,y)$ and $\bfx=(x,y,z)$.  We search for propagating modes along the axis of the microstrip or waveguide having a given angular frequency $\omega>0$.  In particular if $({\cal E},{\cal H})$ denotes the three-dimensional time harmonic electric and magnetic fields in the device, we assume that
\begin{equation}
{\cal E}(\bfx)=\bfE(\hat\bfx)\exp(i\beta z)\mbox{ and }{\cal H}(\bfx)=\bfH(\hat\bfx)\exp(i\beta z)\label{waveguide_form}
\end{equation}
where $\beta$ is the aforementioned propagation constant. Note that if ${\bf\cal U}(\bfx)=\bfU(\hat{\bfx})\exp(i\beta z)$ for some $\bfU=( U_1, U_2, U_3)^T\in \mathbb{R}^3$ then
\[
\nabla \times{\bf\cal U}=\left(\begin{array}{c}{U}_{3,y}-i\beta {U}_2\\-{U}_{3,x}+i\beta {U}_1\\{U}_{2,x}-{U}_{1,y}\end{array}\right)\exp(i\beta z).
\]
So the first order Maxwell system for $({\bf\cal E},{\bf\cal H})$ can be rewritten as
\begin{subequations}\label{beta1}
\begin{eqnarray}
i\omega\epsilon \left(\begin{array}{c}{E}_1\\{E}_2\\{E}_3\end{array}\right)&=&-\left(\begin{array}{c}{H}_{3,y}-i\beta {H}_2\\-{H}_{3,x}+i\beta {H}_1\\{H}_{2,x}-{H}_{1,y}\end{array}\right),\label{beta1_1}\\
i\omega\mu \left(\begin{array}{c}{H}_1\\{H}_2\\{H}_3\end{array}\right)&=&\left(\begin{array}{c}{E}_{3,y}-i\beta {E}_2\\-{E}_{3,x}+i\beta {E}_1\\{E}_{2,x}-{E}_{1,y}\end{array}\right)\label{beta1_2}.
\end{eqnarray}
\end{subequations}
We now set
\[
\hboldE=\left(\begin{array}{c}{E}_1\\{E}_2\end{array}\right) \mbox{ and } \hboldH=\left(\begin{array}{c}{H}_1\\{H}_2\end{array}\right)
\]
and recall the standard scalar and vector curl in the plane given by
\[
\mbox{curl}\,\hat\boldw=w_{2,x}-w_{1,y}\quad\mbox{ and }\quad\mbox{Curl} \,w=\left(\begin{array}{c}w_y\\-w_x\end{array}\right)
\]
where $\hat{\boldw}=(w_1,w_2)^T$.
Using the last sub-equation of \eqref{beta1_2} and the first two equations of \eqref{beta1_1} we have
\[
\Curl\frac{1}{\mu}\curl \hboldE=i\omega\Curl H_3=\omega^2\epsilon\hboldE-\omega\beta \left(\begin{array}{c} H_2\\-H_1\end{array}\right).
\]
Then using the first two equations of \eqref{beta1_2} this can be rewritten as
\begin{equation}
\Curl\frac{1}{\mu}\curl \hboldE-\omega^2\epsilon\hboldE+\frac{\beta^2}{\mu}\hboldE+\frac{i\beta}{\mu}\nabla{E}_3=0.
\label{eq:VD1}
\end{equation}
Next using the first two sub-equations of \eqref{beta1_2} we obtain
\begin{eqnarray*}
{E}_{3,x}=i\beta {E}_1 -i\omega\mu {H}_2,&\quad&
{E}_{3,y}=i\beta {E}_2 +i\omega\mu {H}_1,
\end{eqnarray*}
so that
\[
\div\frac{1}{\mu}\nabla {E}_3=i\beta\div\frac{1}{\mu}\hboldE+i\omega({H}_{1,y}-{H}_{2,x}).
\]
Then using the third sub-equation of \eqref{beta1_1} we obtain 
\begin{equation}
\div\frac{1}{\mu}\nabla {E}_3+\omega^2\epsilon {E}_3-i\beta\div\frac{1}{\mu}\hboldE=0.
\label{eq:VD2}
\end{equation}
Finally, the fact that the 3D field is divergence free gives
\begin{equation}
\div\epsilon\hboldE+i\beta \epsilon {E}_3=0\label{eq:VD3}.
\end{equation}
Equations \eqref{eq:VD1}-\eqref{eq:VD3} constitute equation (2.2) from \cite{VardapetyanDemkowicz:03} allowing for $j=-i$.  They hold in $\Omega$. As mentioned in the introduction, we  assume that the  domain is enclosed in a perfectly conducting (PEC) boundary.  In addition there may be PEC surfaces within the domain representing thin metallic components (if the thickness $t$ of the strip in Fig.~\ref{Fig1} is very small for example). 
On the PEC surface and the transmission line, we have $\nu_1{E}_2-\nu_2{E}_1=0$ where $(\nu_1,\nu_2)^T$ is the unit normal on these edges.  In addition, for the PEC surfaces we know that ${E}_3=0$ so that the above equations are
supplemented by the boundary conditions
\begin{equation}
\bftv\cdot\hboldE=0\mbox{ and } {E}_3=0\mbox{ on }\partial \Omega,\label{eq:VDBC}
\end{equation}
where $\boldsymbol{\tau}=(-\nu_2,\nu_1)^T$.
However equations  \eqref{eq:VD1}-\eqref{eq:VD3}  are redundant.  
The eigenvalue problem that is the main focus of \cite{VardapetyanDemkowicz:03} uses  \eqref{eq:VD1} and \eqref{eq:VD3}
and seeks $\beta\in \mathbb{C}$ and non-trivial $(\boldE,E_3)\in H_0(\curl;\Omega)\times H^1_0(\Omega)$ such that
\begin{subequations}\label{eq:systemVD1}
\begin{align}
\Curl\frac{1}{\mu}\curl \hboldE-\omega^2\epsilon\hboldE+\frac{\beta^2}{\mu}\hboldE+\frac{i\beta}{\mu}\nabla{E}_3&=0
\mbox{ in }\Omega,\label{eq:systemVD1a}\\
\div\epsilon\hboldE+i\beta \epsilon {E}_3&=0 \mbox{ in }\Omega.\label{eq:systemVD1b}
\end{align}
\end{subequations}
An alternative formulation is to use
\eqref{eq:VD1} and \eqref{eq:VD2} so seeking $\beta\in \mathbb{C}$ and non-trivial $(\hboldE,E_3)\in H_0(\curl;\Omega)\times H^1_0(\Omega)$ such that
\begin{subequations}\label{eq:systemVD2}
\begin{align}
\Curl\frac{1}{\mu}\curl \hboldE-\omega^2\epsilon\hboldE+\frac{\beta^2}{\mu}\hboldE+\frac{i\beta}{\mu}\nabla{E}_3&=0
\mbox{ in }\Omega,\\
\div\frac{1}{\mu}\nabla {E}_3+\omega^2\epsilon {E}_3-i\beta\div\frac{1}{\mu}\hboldE&=0\mbox{ in }\Omega.\label{eq:systemVD2b}
\end{align}
\end{subequations}
This formulation will be used to prove the density of the eigenmodes in Section~\ref{density}.

Note that if $(\beta,(\hboldE,E_3))$ is an eigenpair, then, assuming $\epsilon$ and $\mu$ are real,
so are $(\ol{\beta},(\ol{\hboldE},-\ol{E}_3))$
and $(-\beta,(\hboldE,-E_3))$.
Hence the eigenvalues are symmetric about both the real and the imaginary axes.

Under our assumption that $\omega\not =0$, if in addition we are not at a cut-off wave number so that $\beta\not=0$,  equations \eqref{eq:systemVD2} follow from \eqref{eq:systemVD1} as can be seen by taking the divergence of \eqref{eq:systemVD1a} and using the result to replace $\div\epsilon \hboldE$ in \eqref{eq:systemVD1a} to obtain \eqref{eq:systemVD2b}.  Proceeding similarly, when $\beta\not=0$, we see that \eqref{eq:systemVD2} implies \eqref{eq:systemVD1}.  

To obtain a linear eigenvalue problem, it is suggested in \cite{VardapetyanDemkowicz:03} to define $\tilde{E}_3=i\beta E_3$ so that system \eqref{eq:systemVD1} becomes 
\begin{subequations}\label{eq:systemVD1m}
\begin{align}
\Curl\frac{1}{\mu}\curl \hboldE-\omega^2\epsilon\boldE+\frac{\beta^2}{\mu}\hboldE+\frac{1}{\mu}\nabla{\tilde{E}}_3&=0
\mbox{ in }\Omega,\label{eq:systemVD1ma}\\
\div\epsilon\hboldE+ \epsilon {\tilde{E}}_3&=0 \mbox{ in }\Omega,\label{eq:systemVD1mb}
\end{align}
\end{subequations}
together with the boundary conditions \eqref{eq:VDBC}.  The numerical analysis of this problem is the subject of  \cite{VardapetyanDemkowicz:03}, and
this is the version they advocate for numerical purposes.
Under the same change of variables, the alternative formulation (\ref{eq:systemVD2}) becomes: seek $\beta\in \mathbb{C}$ and non-trivial $(\hboldE,\tilde{E}_3)\in H_0(\curl;\Omega)\times H^1_0(\Omega)$ such that
\begin{subequations}\label{eq:systemVD2m}
\begin{align}
\Curl\frac{1}{\mu}\curl \hboldE-\omega^2\epsilon\hboldE+\frac{\beta^2}{\mu}\hboldE+\frac{1}{\mu}\nabla{\tilde{E}}_3&=0
\mbox{ in }\Omega,\\
\div\frac{1}{\mu}\nabla {\tilde{E}}_3+\omega^2\epsilon \tilde{E}_3+\beta^2\div\frac{1}{\mu}\hboldE&=0,\mbox{ in }\Omega,\label{eq:systemVD2bm}
\end{align}
\end{subequations}
together with the boundary conditions \eqref{eq:VDBC}. Now the two approaches 
define the same eigenvalues and eigenfunctions even if $\beta=0$.  

\section{Density of eigenmodes}\label{density}

This section is devoted to analyzing the eigenmodes using \eqref{eq:systemVD2} together with the boundary conditions \eqref{eq:VDBC}. 
We assume that the equation $-\div\mu^{-1}\nabla E_3-\omega^2\epsilon E_3=0$ in $\Omega$, $E_3=0$ on $\partial\Omega$ is uniquely solvable which holds if we exclude those frequencies $\omega$ that support the existence of cut-off wave numbers $\beta=0$.
To perform this analysis we define
\[
X=H_0(\curl,(\div\mu^{-1})^0;\Omega)\times\nabla H^1_0(\Omega)
\]
where, recalling that $\boldsymbol{\tau}$ denotes the unit tangent on $\partial \Omega$,
\[
H_0(\curl,(\div\mu^{-1})^0;\Omega)=\{\bfu\in H({\rm curl};\Omega)\;|\; \boldsymbol{\tau}\cdot \bfu=0\mbox{ on }\partial \Omega\mbox{ and }
\div(\mu^{-1}\bfu)=0\mbox{ in }\Omega\}.
\]
We will use the following scalar product on this space
\[
\spl (\bfu,w),(\bfu',w')\spr_X=\spl \mu^{-1}\curl\boldu,\curl\boldu'\spr
+\spl\epsilon\nabla w,\nabla w'\spr
+\spl\epsilon\boldu,\boldu'\spr,
\]
which simplifies some of the analysis.
This section is devoted to proving the main theorem of this paper:
\begin{theorem}\label{th1} Assume that $\omega$ is not  a cutoff frequency.
Then the spectrum of the eigenvalue problem \eqref{eq:systemVD2m} 
with the boundary conditions \eqref{eq:VDBC} consists of an infinite sequence of eigenvalues, each with finite algebraic multiplicity, which do not accumulate in $\mathbb{C}$  and such that for every $\delta>0$ only finitely many eigenvalues lie outside the sectors $\{z\in \mathbb{C}\;|\;  |\arg z-\pi|<\delta\}$ and $\{z \in\mathbb{C}\;|\; |\arg z+\pi |<\delta\}$.
Finally the closure of the space spanned by  the $\hboldE$ components of the generalized eigenspaces of the problem is dense in $H_0(\curl;\Omega)$. 
\end{theorem}
\begin{remark}
The proof of Theorem~\ref{th1} relies on the fact that certain operators are of finite order.  To define this recall from \cite{Halla:22MSteklovStab} that for Hilbert spaces  ${\cal X},$ and ${\cal Y}$ a compact operator $K\in L({\cal X}, {\cal Y})$ is in Schatten class $K_p({\cal X},{\cal Y})$ of order $p \in (0,\infty)$, if the sequence of singular values $s_n(K),$ $n\in \mathbb{N}$, of $K$ is $\ell^p(\mathbb{N})$ summable. Then, a compact operator
$K \in L({\cal X},{\cal Y})$ is said to be of finite order if there exists $p \in (0,\infty)$ such that $K\in K_p({\cal X},{\cal Y})$.
\end{remark}

\begin{proof} We start by building build the Schur complement with respect to $E_3$ and obtain
\begin{subequations}\label{eq:evp-Schur}
\begin{align}
\Curl\mu^{-1}\curl \hboldE - \omega^2\epsilon\hboldE+\beta^2\mu^{-1}\hboldE\label{Schur1}
+\beta^2\mu^{-1}\nabla S \div\mu^{-1}\hboldE &=0,\\
\tv\cdot\hboldE&=0,\label{Schur2}
\end{align}
\end{subequations}
where formally $S:=(-\div\mu^{-1}\nabla-\omega^2\epsilon)^{-1}\in L(H^{-1}(\Omega),H^1_0(\Omega))$.
Next we use the orthogonal decomposition $H_0(\curl;\Omega)=H_0(\curl,(\div\mu^{-1})^0;\Omega)\oplus^\bot \nabla H^1_0(\Omega)$ and write $\hboldE=\boldu+\nabla w$.
Using this expansion, we obtain the variational equation for $(\bfu,w)\in X$
\begin{align*}
0&=\spl \mu^{-1}\curl\boldu,\curl\boldu'\spr-\omega^2\spl\epsilon\boldu,\boldu'\spr+\beta^2\spl\mu^{-1}\boldu,\boldu'\spr\\
&-\omega^2\big(\spl\epsilon\boldu,\nabla w'\spr+\spl\epsilon\nabla w,\boldu'\spr\big)\\
&-\omega^2 \spl\epsilon\nabla w,\nabla w'\spr+\beta^2 \spl\mu^{-1}\nabla w,\nabla w'\spr
+\beta^2 \spl \mu^{-1}\nabla S\div\mu^{-1}\nabla w,\nabla w'\spr.
\end{align*}
for all $(\bfu',w')\in X$.
We compute, using the definition of $S$
\begin{align*}
\spl \mu^{-1}\nabla S\div\mu^{-1}\nabla w,\nabla w'\spr
&=\spl \mu^{-1}\nabla S(\div\mu^{-1}\nabla w+\omega^2\epsilon w),\nabla w'\spr
-\omega^2\spl \mu^{-1}\nabla S \epsilon w),\nabla w'\spr\\
&=-\spl \mu^{-1}\nabla w,\nabla w'\spr
-\omega^2\spl \mu^{-1}\nabla S \epsilon w),\nabla w'\spr
\end{align*}
Then using integration by parts and the definition of $S$ again:
\begin{align*}
\spl \mu^{-1}\nabla S\div\mu^{-1}\nabla w,\nabla w'\spr
&=-\spl \mu^{-1}\nabla w,\nabla w'\spr
+\omega^2\spl \epsilon w,S\div\mu^{-1}\nabla w'\spr\\
&=-\spl \mu^{-1}\nabla w,\nabla w'\spr
+\omega^2\spl \epsilon w,S(\div\mu^{-1}\nabla w'+\omega^2\epsilon w')\spr
-\omega^4\spl \epsilon w,S \epsilon w')\spr\\
&=-\spl \mu^{-1}\nabla w,\nabla w'\spr
-\omega^2\spl \epsilon w, w'\spr
-\omega^4\spl \epsilon w,S \epsilon w'\spr.
\end{align*}
Hence $(\bfu,w)\in X$
\begin{align*}
&\spl \mu^{-1}\curl\boldu,\curl\boldu'\spr-\omega^2\spl\epsilon\boldu,\boldu'\spr+\beta^2\spl\mu^{-1}\boldu,\boldu'\spr\\
&-\omega^2\big(\spl\epsilon\boldu,\nabla w'\spr+\spl\epsilon\nabla w,\boldu'\spr\big)\\
&-\omega^2 \spl\epsilon\nabla w,\nabla w'\spr
-\beta^2 \omega^2\big(\spl \epsilon w, w'\spr
+\omega^2\spl S\epsilon w, \epsilon w'\spr\big)=0
\end{align*}
for all $(\bfu',w')\in X$.
We switch the sign of the test function $w'\to-w'$ and obtain
\begin{align}
(A-T-\lambda K)(\boldu,w)=0\label{opeq}
\end{align}
with $\lambda:=-\beta^2$ and
\begin{align*}
\spl A(\boldu,w),(\boldu',w')\spr&:=
\spl \mu^{-1}\curl\boldu,\curl\boldu'\spr
+\spl\epsilon\nabla w,\nabla w'\spr
+\spl\epsilon\boldu,\boldu'\spr,\\
\spl T(\boldu,w),(\boldu',w')\spr&:=
\omega^2\big(\spl\epsilon\boldu,\boldu'\spr-\spl\epsilon\boldu,\nabla w'\spr+\spl\epsilon\nabla w,\boldu'\spr\big)
+\spl\epsilon\boldu,\boldu'\spr,\\
\spl K(\boldu,w),(\boldu',w')\spr&:=
\spl\mu^{-1}\boldu,\boldu'\spr
-\omega^2\big(\spl \epsilon w, w'\spr
+\omega^2\spl S\epsilon w, \epsilon w'\spr\big).
\end{align*}
In view of the definition of the scalar product $\spl\cdot,\cdot\spr_X$ we see that $A=I$.

Both $T$ and $K$ are compact and of finite order.
In particular, consider the embedding operator $E\in L(H_0(\curl,(\div\mu^{-1})^0;\Omega),\boldL^2(\Omega))$, the multiplication operator $M_\epsilon\in L(\boldL^2(\Omega))$ and the gradient operator $G\in L(H^1_0(\Omega),\boldL^2(\Omega))$, $Gw:=\nabla$.
Thence $T=\omega^2( E^*M_\epsilon E- G^*M_\epsilon E+ E^*M_\epsilon G) + E^* M_\epsilon E$.
Due to \cite{Ciarlet:20} there exists $s>0$ such that $H_0(\curl,(\div\mu^{-1})^0;\Omega)$ embeds continuously into $\boldH^s(\Omega)$.
Thus it follows with \cite[Thm.~2.2, Lem.~2.3]{Halla:22MSteklovStab} that $E$ and hence $T$ is of finite order.
With the same approach it can be seen that $K$ is of finite order too.
In addition $K$ is selfadjoint.

The operator $K$ is also injective.
Indeed let $(\boldu,w)\in\ker K$.
It easily follows that $\boldu=0$ and $w+\omega^2S\epsilon w=0$.
Multiplying with $S^{-1}$ we obtain
\begin{align*}
0=-\div\mu^{-1}\nabla w-\omega^2\epsilon w+\omega^2\epsilon w=-\div\mu^{-1}\nabla w,
\end{align*}
and thus $w=0$.
Now we can simply apply \cite[Theorem~2.1]{Halla:22MSteklovStab}.
It remains to note that for an eigenpair $(\lambda,(\boldu,w))$ of \eqref{opeq}, $(\pm\sqrt{-\lambda},\boldu+\nabla w)$ are eigenpairs of \eqref{eq:evp-Schur} and vice-versa.

\end{proof}
Note that here we do not have to deal with the involved concept of double completeness, because the eigenfunctions for $\pm\beta$ are identical.
Furthermore, for topologically nontrivial domains the operator $K$ remains injective.

Note also that we could use this formulation also for the computation of modes.
We can enforce the constraint $\div\mu^{-1}\boldu=0$ weakly by the introduction of a scalar auxiliary variable.
The implementation of $S$ requires also the introduction of a scalar auxiliary variable.
Thus we would obtain a problem in the space $H_0(\curl;\Omega)\times H^1_0(\Omega)\times H^1_0(\Omega)$, which is more expensive than 
using the space $H_0(\curl;\Omega)\times H^1_0(\Omega)$ which is used for \eqref{eq:systemVD1}.

\section{A transparent boundary condition}\label{DtN}
Before we start the main derivation of this section we prove the following lemma.
\begin{lemma}\label{lem:orthog}
Assume that the equation $-\div\mu^{-1}\nabla E_3-\omega^2\epsilon E_3=0$ in $\Omega$, $E_3=0$ on $\partial\Omega$ is uniquely solvable.  In addition, suppose $(\beta_j,(\hboldE_j,E_{3,j}))$ and $(\beta_k,(\hboldE_k,E_{3,k}))$ are both eigenpairs of (\ref{eq:systemVD2m}).  Assume $\beta_j^2\not=\overline{\beta}_k^2$, $\beta_j\not=0$, $\beta_k\not=0$, then
\begin{align*}
\spl \mu^{-1}\curl\hboldE_j,\curl\hboldE_k\spr-\omega^2\spl\epsilon\hboldE_j,\hboldE_k\spr=0.
\end{align*}
\end{lemma}

\begin{proof}
Using the Schur complement derived in the proof of the previous theorem (see (\ref{Schur1})-(\ref{Schur2})) we see that 
 the variational formulation of (\ref{eq:systemVD2m}) is to find $\beta$ and $\hboldE\in  H_0(\curl;\Omega)$ with $\hboldE\not=0$ such that
\begin{equation}
\langle \mu^{-1}\curl\hboldE,\curl \xi\rangle-\omega^2\langle \epsilon \hboldE,\xi\rangle=-\beta^2\left(\langle \mu^{-1}\hboldE,\xi\rangle-\langle S\left(\div\frac{1}{\mu}\hboldE\right),\div \frac{1}{\mu}\xi\rangle\right)\label{eq:junk}
\end{equation}
for all $\xi\in  H_0(\curl;\Omega)$.

Equation (\ref{eq:junk}) is of the kind $A\hboldE=\beta^2B\hboldE$ with self-adjoint operators $A,B$.
Hence
$\spl A\hboldE_k,\hboldE_j \spr = \beta_k^2 \spl B\hboldE_k,\hboldE_j \spr$
and also
$\spl A\hboldE_k,\hboldE_j \spr = \spl \hboldE_k,A\hboldE_j \spr =  
\ol{\beta_j}^2 \spl \hboldE_k,B\hboldE_j \spr
=\ol{\beta_j}^2 \spl B\hboldE_k,\hboldE_j \spr$.
Thus $(\beta_k^{-2}-\ol{\beta_j}^{-2}) \spl A\hboldE_k,\hboldE_j \spr = 0$.
If $(\beta_k^{-2}-\ol{\beta_j}^{-2})$ does not vanish, $\spl A\hboldE_k,\hboldE_j \spr$ has to vanish.
\end{proof}
\begin{figure}
\centering
\resizebox{0.5\textwidth}{!}{\includegraphics{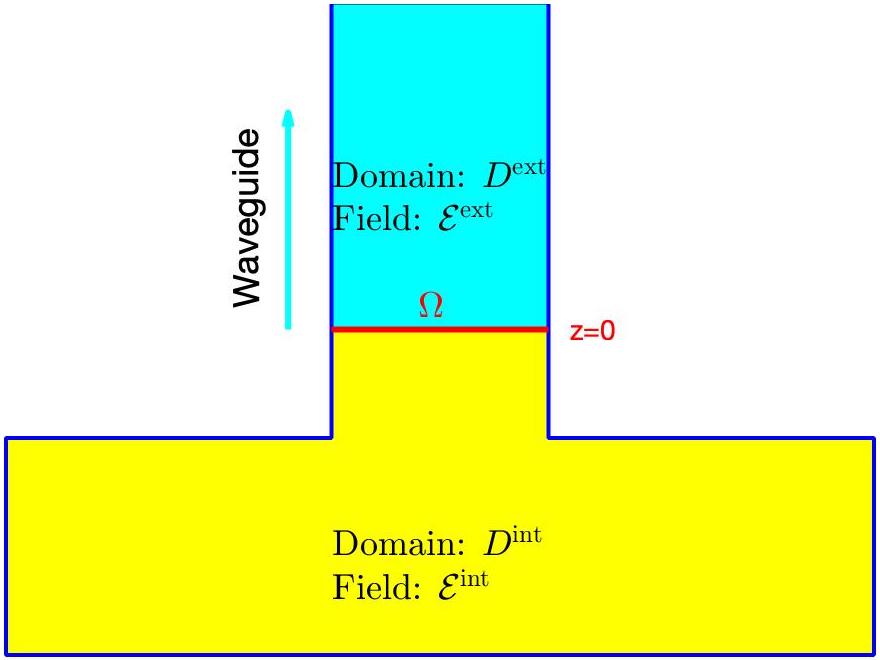}}
\caption{Cartoon of a terminated wave guide.  The semi-infinite waveguide $D^{\rm ext}$ (blue region) is terminated by the bounded domain $D^{\rm int}$ (yellow region).
If finite elements are used to discretize the termination region $D^{\rm int}$, we can provide a boundary condition on the artificial boundary $\Omega$
by using a Dirichlet-to-Neumann map constructed from the modes propagating in the waveguide.}
\label{fig2}
\end{figure}
Now we can show how the eigen-modes can be used to provide a termination condition for waveguide calculations.
When simulating waveguides by finite elements  it is  necessary to truncate the unbounded computational domain.  One way to do this is using a transparent boundary condition
that can be computed using the modes we have analyzed in the previous section~\cite{Kim:17}.  A cartoon example is shown in Fig.~\ref{fig2} where a vertical waveguide denoted $D^{\mathrm{ext}}$ is terminated by a bounded domain $D^{int}$.   We assume a source ${\cal J}$ is located in the bounded domain $D^{\mathrm{\mathrm{int}}}$ (having compact support) and let $(\epsilon,\mu)$ denote the electromagnetic parameters of the three dimensional domain.
In $D^{\mathrm{ext}}$, the functions $\mu$ and $\epsilon$ are assumed to satisfy the assumptions at the start of Section~\ref{deriv} so that $D^{\mathrm{ext}}$ is a waveguide and we can use the 
modes analyzed in the previous section to represent the solution on $\Omega$ (due to their verified density).  In $D^{\rm int}$ the electromagnetic parameters are positive real piecewise $W^{1,\infty}$ functions.  Using finite element elements the termination region $D^{int}$  can be discretized using finite elements.  Further up the waveguide we need to terminate the computational domain by cutting the waveguide at $\Omega$ which for simplicity we assume is in the plane $z=0$. Then the total electric field ${\cal E}^{\mathrm{int}}$ (now a general vector function not of the form (\ref{waveguide_form})) satisfies
\[
\nabla_3\times \mu^{-1}\nabla_3\times {\cal E}^{\mathrm{int}}-\omega^2\epsilon{\cal E}^{\mathrm{int}}=i\omega {\cal J}\mbox{ in } D^{\mathrm{int}}\subset \mathbb{R}^3,
\]
together with the PEC boundary condition on the external boundary of $D^{\mathrm{int}}$. Here $\nabla_3\times$ is the full three dimensional curl.
Multiplying by a smooth test vector function ${\cal U}'$ and integrating by parts, we obtain the standard variational problem
\[
\int_{D^{\mathrm{int}}}\mu^{-1}\nabla_3\times {\cal E}^{\mathrm{int}}\cdot\nabla_3\times \overline{\cal U}'-\omega^2\epsilon{\cal E}^{\mathrm{int}}\cdot\overline{\cal U}'\,dV+\int_\Omega \bfe_3\times \mu^{-1}\nabla_3\times {\cal E}^{\mathrm{int}}\cdot\overline{\cal U}'\,dA=
\int_{D^{\mathrm{int}}}i\omega {\cal J}\cdot\overline{\cal U}'\,dV,
\]
where $\bfe_3=(0,0,1)^T$.
By continuity of the tangential trace of the electric and magnetic fields across $\Omega$,
\[
\int_\Omega \bfe_3\times \mu^{-1}\nabla_3\times {\cal E}^{\mathrm{int}}\cdot\overline{\cal U}'\,dA=\int_\Omega \bfe_3\times \mu^{-1}\nabla_3\times {\cal E}^{\mathrm{ext}}\cdot\overline{\cal U}'\,dA.
\]
 In $D^{\mathrm{ext}}$ we set
\begin{align*}
\calE^\mathrm{ext}(\boldx,z)=\left( \begin{array}{c}\hboldE^\mathrm{ext}( \boldx,z)\\E_3^\mathrm{ext}( \boldx,z) \end{array}\right).
\end{align*} 
In the domain $D^{\rm ext}\Omega\times (,\infty)$, the field $\calE^\mathrm{ext}$ satisfies Maxwell's equations and is outgoing.  In addition $\bfe_3\times(\calE^\mathrm{ext}\times \bfe_3)=\bfe_3\times(\calE^\mathrm{int}\times\bfe_3)$ at $z=0$.  By the density of the modes on $\Omega$ we
can approximate 
\begin{equation}
\bfe_3\times(\calE^\mathrm{int}(\bfx,0)\times\bfe_3)\approx \sum_{j=1}^N a_j\hboldE_j(\bfx)\mbox{ in } H_0({\rm curl};\Omega)
\label{exp}
\end{equation}
where we can achieve any desired accuracy by taking $N$ large enough.  The equality of the traces above and below $\Omega$ implies that
\[
\calE^\mathrm{ext}(\bfx,z)\approx \sum_{j=1}^N a_j\left(\begin{array}{c}\hboldE_j(\bfx)\\E_{3,j}(\bfx)\end{array}\right)\exp(i\beta_j z)
\]
where $\beta_j$ is chosen such that $\beta_j>0$ if $\beta_j$ is real, and $\Im(\beta_j)>0$ if $\beta_j$ has a non-zero imaginary part in order to obtain an outgoing or evanescent mode.
Then, as we have seen, we encounter
\begin{align*}
\int_{\Omega}\bfe_z\times\mu^{-1}\nabla_{3}\times\calE^\mathrm{ext}\cdot\overline{\calU}'\,dA
&=\left\spl \bfe_z\times\mu^{-1} \left(\begin{array}{c}E^\mathrm{ext}_{3,y}-E^\mathrm{ext}_{2,z}\\-E^\mathrm{ext}_{3,x} + E^\mathrm{\mathrm{int}}_{1,z}\\E^\mathrm{ext}_{2,x}-E^\mathrm{ext}_{1,y}\end{array}\right)
, \calU'\right\spr\\
&=\left\spl \mu^{-1} \left(\begin{array}{c}E^\mathrm{ext}_{3,x}- E^\mathrm{ext}_{1,z}\\E^\mathrm{ext}_{3,y}-E^\mathrm{ext}_{2,z}\\0\end{array}\right), \calU'\right\spr\\ &=\left\spl\mu^{-1} \bpm \nabla E^{\rm ext}_3-\frac{\partial}{\partial z}\hboldE^{\rm ext}\\0\epm,{\cal U}'\right\spr.
\end{align*}
To construct a transparent boundary condition by means of a Dirichlet-to-Neumann operator we need to derive a meaningful expression for the right hand side above.
For a mode of the cross sectional problem $(\hboldE,E_3)$, with propagation constant $\beta$ assuming that $\beta\not=0$,  we note that the corresponding wave-guide mode is
$(\hboldE,E_3)\exp(i\beta z)$ so for this mode we compute
\begin{align*}
\mu^{-1} \left(\begin{array}{c}E_{3,x}- E_{1,z}\\E_{3,y}-E_{2,z}\\0\end{array}\right)
&=\mu^{-1} \bpm \nabla E_3-i\beta\hboldE \\ 0 \epm\exp(i\beta z)\\
&=\mu^{-1} \bpm \frac{\mu}{i\beta}\left(-\Curl\frac{1}{\mu}\curl \hboldE+\omega^2\epsilon\hboldE-\frac{\beta^2}{\mu}\hboldE)\right)-i\beta\hboldE \\ 0 \epm\exp(i\beta z)\\
&= i\beta^{-1} \bpm  \Curl\mu^{-1}\curl\hboldE-\omega^2\epsilon\hboldE \\ 0 \epm\exp(i\beta z).
\end{align*}
Note that $\calU'$  has a vanishing tangential trace on the waveguide boundary, and the tangential trace $\hboldE^\mathrm{int}$ of ${\cal E}^{\mathrm{int}}$ on the interface $\Omega$ has a vanishing tangential trace on $\partial\Omega$.
Assuming  that $\beta=0$ is not an eigenvalue, we can normalize the eigenfunctions such that
\begin{align*}
\spl \mu^{-1}\curl\hboldE_j,\curl\hboldE_j\spr-\omega^2\spl\epsilon\hboldE_j,\hboldE_j\spr
=1.
\end{align*}
We then need to compute an expansion of the field in terms of modes on $\Omega$ (i.e. find the coefficients $a_n$ in the expansion (\ref{exp}). To see how this can be done, suppose that $\beta_j^2$, $j=1,\cdots,\infty$ are real and the eigenvalues are simple. Then the propagation constants are either real or purely imaginary (i.e. no real part), so using Lemma~\ref{lem:orthog} the projection of $\hboldE^{\rm int}$  onto $\hboldE_j$ is given by
\begin{align*}
\big(\spl \mu^{-1}\curl\hboldE^{\rm int},\curl\hboldE_j\spr-\omega^2\spl\epsilon\hboldE^{\rm int},\hboldE_j\spr\big)\hboldE_j.
\end{align*}
Thus, in this simple case, to construct an (approximate) DtN-operator, we have the following steps
\begin{enumerate}
 \item Choose $N$ and compute the modes $(\beta_j,\boldE_j)$, ($j=1,\dots,N$).
 \item Compute the boundary term
 \begin{align*}
 \sum_j -i\beta_j^{-1} \big(\spl \mu^{-1}\curl\hboldE^{\rm int},\curl\hboldE_j\spr-\omega^2\spl\epsilon\hboldE^{\rm int},\hboldE_j\spr\big)
 \big(\spl \mu^{-1}\curl\hboldE_j,\curl{\cal U'}\spr-\omega^2\spl\epsilon\hboldE_j,{\cal U}'\spr\big).
 \end{align*}
\end{enumerate}
If any of the propagation constants $\{\beta_j^2\}_{j=1}^{\infty}$ are complex or repeated eigenvalues occur, then the projection must be computed more carefully. Fixing $j$, let ${\cal K}_j$ denote the set of all modes with index $k$ such that $\beta_j^2=\overline{\beta}_k^2$.  Then the modes may not be orthogonal, and to compute the expansion of $\hboldE^{\rm int}$ of ${\cal K}_j$ would require solving a small linear system.

\section*{Acknowledgements}
The first author of this work was supported by DFG project 468728622 and DFG SFB 1456 project 432680300.
\bibliographystyle{amsplain}
\bibliography{short_biblio}

\providecommand{\bysame}{\leavevmode\hbox to3em{\hrulefill}\thinspace}
\providecommand{\MR}{\relax\ifhmode\unskip\space\fi MR }
\providecommand{\MRhref}[2]{%
  \href{http://www.ams.org/mathscinet-getitem?mr=#1}{#2}
}
\providecommand{\href}[2]{#2}
\begin{thebibliography}{10}

\bibitem{Bamberger:90}
A.~Bamberger and A.S. Bonnet, \emph{Mathematical analysis of the guided modes
  of an optical fiber}, SIAM J. Math. Anal. \textbf{21} (1990), 1487--1510.

\bibitem{Ciarlet:20}
Patrick Ciarlet, \emph{{On the approximation of electromagnetic fields by edge
  finite elements. III: Sensitivity to coefficients}}, {SIAM J. Math. Anal.}
  \textbf{52} (2020), no.~3, 3004--3038.

\bibitem{GohbergKrein:69}
I.~C. Gohberg and M.~G. Kre\u{\i}n, \emph{Introduction to the theory of linear
  nonselfadjoint operators}, Translated from the Russian by A. Feinstein.
  Translations of Mathematical Monographs, Vol. 18, American Mathematical
  Society, Providence, R.I., 1969.

\bibitem{Halla:22MSteklovStab}
Martin Halla, \emph{On the existence and stability of modified {M}axwell
  {S}teklov eigenvalues}, 2022, https://arxiv.org/abs/2207.06498.

\bibitem{JolyPoirier:95}
Patrick Joly and Christine Poirier, \emph{Mathematical analysis of
  electromagnetic open waveguides}, RAIRO, Mod{\'e}lisation Math. Anal.
  Num{\'e}r. \textbf{29} (1995), no.~5, 505--575.

\bibitem{JolyPoirier:99}
\bysame, \emph{A numerical method for the computation of electromagnetic modes
  in optical fibres}, Math. Methods Appl. Sci. \textbf{22} (1999), no.~5,
  389--447.

\bibitem{JolyPoirierRobertsTrouve:96}
Patrick Joly, Christine Poirier, J.~E. Roberts, and P.~Trouve, \emph{A new
  nonconforming finite element method for the computation of electromagnetic
  guided waves. {I}: {Mathematical} analysis}, SIAM J. Numer. Anal. \textbf{33}
  (1996), no.~4, 1494--1525.

\bibitem{Kim:17}
Seungil Kim, \emph{Analysis of the non-reflecting boundary condition for the
  time-harmonic electromagnetic wave propagation in waveguides}, Journal of
  Mathematical Analysis and Applications \textbf{453} (2017), no.~1, 82--103.

\bibitem{Lumerical}
Ansys Lumerical, \emph{Microstrip transmission line}, \\
  \url{https://support.lumerical.com/hc/en-us/articles/360042051394-Microstrip-transmission-line},
  Accessed: 2/21/2022.

\bibitem{Markus:88}
A.~S. Markus, \emph{Introduction to the spectral theory of polynomial operator
  pencils}, Translations of Mathematical Monographs, vol.~71, American
  Mathematical Society, Providence, RI, 1988. \MR{971506}

\bibitem{Pozar:11}
D.M. Pozar, \emph{Microwave engineering}, 4 ed., John Wiley \& Sons, inc, 2011.

\bibitem{SheSmi13a}
Yury Shestopalov and Yury Smirnov, \emph{Eigenwaves in waveguides with
  dielectric inclusions: completeness}, Applicable Analysis \textbf{93} (2013),
  1824--1845.

\bibitem{SheSmi13}
\bysame, \emph{Eigenwaves in waveguides with dielectric inclusions: spectrum},
  Applicable Analysis \textbf{93} (2013), 408--427.

\bibitem{Dassault}
Dassault Simulia, \emph{Microstrip transmission line simulation}, \\
  \url{https://www.3ds.com/products-services/simulia/resources/microstrip-transmission-line/},
  Accessed: 2/21/2022.

\bibitem{VardapetyanDemkowicz:03}
L.~Vardapetyan and L.~Demkowicz, \emph{Full-wave analysis of dielectric
  waveguides at a given frequency}, Math. Comput. \textbf{72} (2003), no.~241,
  105--129.

\bibitem{VardapetyanDemkowiczNeikirk:03}
L.~Vardapetyan, L.~Demkowicz, and D.~Neikirk, \emph{{{\(hp\)}}-vector finite
  element method for eigenmode analysis of waveguides}, Comput. Methods Appl.
  Mech. Eng. \textbf{192} (2003), no.~1-2, 185--201.

\end{thebibliography}

\end{document}